\documentclass{amsart}

\usepackage{amsmath,amssymb,amscd,amsfonts}

\newtheorem{theorem}{Theorem}
\newcommand{\bt}{\begin{theorem}}
\newcommand{\et}{\end{theorem}}
\newtheorem{lemma}{Lemma}
\newcommand{\bl}{\begin{lemma}}
\newcommand{\el}{\end{lemma}}
\newtheorem{corollary}{Corollary}
\newcommand{\bc}{\begin{corollary}}
\newcommand{\ec}{\end{corollary}}
\newtheorem{problem}{Problem}
\newcommand{\bprob}{\begin{problem}}
\newcommand{\eprob}{\end{problem}}

\newtheorem{example}{Example}
\newcommand{\bex}{\begin{example}}
\newcommand{\eex}{\end{example}}

\newcommand{\beq}{\begin{equation}}
\newcommand{\eeq}{\end{equation}}
\newcommand{\benum}{\begin{enumerate}}
\newcommand{\eenum}{\end{enumerate}}
\newcommand{\N}{\ensuremath{ \mathbf N }}
\newcommand{\Z}{\ensuremath{\mathbf Z}}
\newcommand{\Q}{\ensuremath{\mathbf Q}}

\newcommand{\mcf}{\ensuremath{ \mathcal F}}

\newcommand{\mba}{\ensuremath{ \mathbf a}}
\newcommand{\mbb}{\ensuremath{ \mathbf b}}

\newcommand{\mbg}{\ensuremath{ \mathbf g}}
\newcommand{\mbh}{\ensuremath{ \mathbf h}}

\newcommand{\mbw}{\ensuremath{ \mathbf w}}
\newcommand{\mbx}{\ensuremath{ \mathbf x}}

\newcommand{\mbz}{\ensuremath{ \mathbf z}}

\newcommand{\Rn}{\ensuremath{ \mathbf{R}^n }}
\newcommand{\bsmallmat}{\left(\begin{smallmatrix}}
\newcommand{\esmallmat}{\end{smallmatrix}\right)}

\DeclareMathOperator{\kernel}{\text{kernel}}

\DeclareMathOperator{\Odd}{\mathbf{Odd}}

\newcommand{\bmat}{\left(\begin{matrix}}
\newcommand{\emat}{\end{matrix}\right)}

\DeclareMathOperator{\rank}{\text{rank}}

\DeclareMathOperator{\vectorsmallzn+1}{\left( \begin{smallmatrix} z_1 \\ \vdots \\ z_n\\ z_{n+1}  \end{smallmatrix}\right)}
\DeclareMathOperator{\vectorsmallz'n+1}{\left( \begin{smallmatrix} z'_1 \\ \vdots \\ z'_n\\ z'_{n+1}  \end{smallmatrix}\right)}

\DeclareMathOperator{\vectorsmallwn+1}{\left( \begin{smallmatrix} w_1 \\ \vdots \\ w_n \\ w_{n+1} \end{smallmatrix}\right)}

\DeclareMathOperator{\vectorsmallxn+1}{\left( \begin{smallmatrix} x_1 \\ \vdots \\ x_n \\ x_{n+1} \end{smallmatrix}\right)}

\DeclareMathOperator{\vectorsmallzn+1}{\left( \begin{smallmatrix} z_1 \\ \vdots \\ z_n\\ z_{n+1}  \end{smallmatrix}\right)}

\newcommand{\mbo}{\ensuremath{\mathbf 0}}

\DeclareMathOperator{\vectorsmallhn}{\left( \begin{smallmatrix} h_1 \\ \vdots \\ h_n \end{smallmatrix}\right)}

\title[Increasing-decreasing patterns of arithmetic functions]{Increasing-decreasing patterns in the iteration of an arithmetic function}
\author{Melvyn B. Nathanson}
\address{Lehman College (CUNY), Bronx, New York 10468}
\email{melvyn.nathanson@lehman.cuny.edu}



\thanks{Supported in part by PSC-CUNY Research Award Program grant 66197-00 54.} 

\date{\today}

\begin{document}

\maketitle

\begin{abstract} 
Let $\Omega$ be a set of positive integers and let $f:\Omega \rightarrow \Omega$
 be an arithmetic function.  Let $V = (v_i)_{i=1}^n$ be a finite sequence of positive integers.  
An integer $m \in \Omega$ has \textit{increasing-decreasing pattern} $V$ with respect to $f$ if,  
for all odd integers $i \in \{1,\ldots, n\}$,  
\[
f^{v_1+ \cdots + v_{i-1}}(m) < f^{v_1+ \cdots + v_{i-1}+1}(m) < \cdots < f^{v_1+ \cdots + v_{i-1}+v_{i}}(m)
\]
and, for all even  integers $i \in \{2,\ldots, n\}$, 
\[
f^{v_1+ \cdots + v_{i-1}}(m) > f^{v_1+ \cdots  +v_{i-1}+1}(m) > \cdots > f^{v_1+ \cdots  +v_{i-1}+v_i}(m).
\]
The arithmetic function $f$ is \textit{wildly increasing-decreasing} if,  
for every finite sequence $V$ of positive integers, there exists an integer $m  \in \Omega$ 
such that $m$ has increasing-decreasing pattern $V$ with respect to $f$.  
This paper gives a proof that the Syracuse function 
is wildly increasing-decreasing.
\end{abstract}

\section{Iterations and patterns}                    \label{Syracuse:section:patterns} 

An arithmetic function is a function whose domain is a subset $\Omega$ of 
the set $\N = \{1, 2, 3, \ldots \}$ of positive integers.   
Let $f: \Omega \rightarrow \Omega$ be an arithmetic function 
and, for all $j \in \N$, let $f^j:\Omega \rightarrow \Omega$  be the $j$th iterate of $f$.   
The function $f^0:\Omega \rightarrow \Omega$  is the identity function.  
The \emph{trajectory}\index{trajectory} of $m \in \Omega$ 
is the sequence of positive integers $\left( f^j(m) \right)_{j=0}^{\infty}$. 
The pair $(\Omega,f)$ is a discrete dynamical system.

 Let 
\[
\Omega^{\text{fix}} = \{m \in \Omega: f(m) = m \}
\]
 be the set of fixed points of $f$ 
and let 
\[
\Omega^{\text{per}} = \{m \in \Omega: f^k(m) = m \text{ for some positive integer $k$} \} 
\]
be the set of periodic points of $f$.   The \emph{period}\index{period} of $m \in \Omega^{\text{per}}$ 
is the smallest positive integer $k$ such that $f^k(m) = m$.
The fixed points are the points of period one.

The  trajectory $\left( f^j(m) \right)_{j=0}^{\infty}$ is  \emph{eventually constant} 
if   $f^k(m) \in \Omega^{\text{fix}}$ for some $k \in \N$. 
The  trajectory $\left( f^j(m) \right)_{j=0}^{\infty}$ is  \emph{eventually periodic} 
if $f^k(m) \in \Omega^{\text{per}}$ for some $k \in \N$.   
The trajectory is  \emph{bounded} if and only if $f^k(m) \in \Omega^{\text{per}}$ for some $k \in \N$.  
The trajectory is  \emph{unbounded} if and only if $\lim_{k\rightarrow \infty} f^k(m) = \infty$.  
If $m \in \Omega$ and $m \notin \Omega^{\text{fix}}$, 
then either $f(m) > m$ (that is, $f$ increases $m$) 
or $f(m) < m$  (that is, $f$ decreases $m$) .  
We are interested in the pattern of changes (increases and decreases) in the trajectory 
$\left( f^j(m) \right)_{j=0}^{\infty}$.

Here are some examples.  
Let $\Omega$ be the set of odd positive integers and let $\ell\in \Omega$.   
For all $m \in \Omega$, the integer $\ell m+1$ is even.  
We define the arithmetic function 
$S_{\ell}: \Omega \rightarrow \Omega$  by 
\[
S_{\ell}(m) = \frac{\ell m+1}{2^e}
\]
where $2^e$ is the highest power of 2 that divides $\ell m+1$. 
The function $S_{\ell}$ has a fixed point if and only if $\ell = 2^k-1$ 
for some positive integer $k$, and in this case the unique fixed point is $m=1$.
The function $S_5$ has no fixed point, but it does have periodic points.  
For example, $S_5(1) = 3$ and $S_5(3) = 1$. 

The \emph{Syracuse function}\index{Syracuse function} 
(also called the \emph{Collatz function}) is the arithmetic function  
\beq           \label{Syracuse:function}
S(m) = S_3(m) = \frac{3m+1}{2^e}
\eeq
where $2^e$ is the highest power of 2 that divides $3m+1$.  
The unique fixed point of this function is $m=1$.  There is a large 
literature on the Syracuse function (cf. Eliahou~\cite{elia93}, Everett~\cite{ever77}, Lagarias~\cite{laga85,laga10}, Tao~\cite{tao22}, Terras~\cite{terr76}). 

Let $\Omega = \N$.    
The \emph{sum of the proper divisors function} is the  arithmetic function  
 $s:\N \rightarrow \N$ defined by 
\[
s(m) = \sum_{\substack{d|m \\ 1 \leq d < m }} d. 
\]
The trajectory of the positive integer $m$, that is, $\left(s^j(m)\right)_{j=1}^{\infty}$, 
 is called 
an \emph{aliquot sequence}\index{aliquot sequence}.  
The fixed points of $s$ are the \emph{perfect numbers}\index{perfect number} $6, 28, 496, \ldots$.  
If $m$ is a point of period 2 and if $n = s(m)$, then $m \neq n$ and $s(n) = m$.  
The pair of integers $(m,n)$ 
is called an \emph{amicable pair}\index{amicable pair}.   For example, $(220, 284)$ 
and $(1184, 1210)$ are amicable pairs.  Integers whose trajectories are periodic under $s$ are 
called \emph{sociable numbers}\index{sociable numbers}.  There exist sociable numbers of period 4 
(such as 1264460), 
but no sociable number of period 3 is known.  
The behavior of aliquot sequences is poorly understood (cf. Guy~\cite{guy04}, 
P. Erd\H os, A. Granville, C. Pomerance, and C. Spiro~\cite{egps90}, 
Pollack-Pomerance~\cite{poll-pome16}). 

Let $\Omega$ be a set of positive integers and let $f: \Omega \rightarrow \Omega$ be an arithmetic function.
Let $V = (v_i)_{i=1}^n$ be a finite sequence of positive integers.  
We say that an integer $m$ in $\Omega$ has \emph{increasing-decreasing pattern}\index{pattern} 
$V$ with respect to $f$ if 
\[
m < f(m) < f^2(m) < \cdots < f^{v_1}(m) 
\]
\[
f^{v_1}(m) > f^{v_1+1}(m) > \cdots > f^{v_1+v_2}(m)
\]
\[
f^{v_1+v_2}(m) < f^{v_1+v_2+1}(m) < \cdots < f^{v_1+v_2+v_3}(m)
\]
and, in general, if $i$ is odd, then 
\beq                                                                 \label{Syracuse:increasing} 
f^{v_1+ \cdots + v_{i-1}}(m) < f^{v_1+ \cdots + v_{i-1}+1}(m) < \cdots < f^{v_1+ \cdots + v_{i-1}+v_{i}}(m)
\eeq
and if $i$ is even, then 
\beq                                                                 \label{Syracuse:decreasing} 
f^{v_1+ \cdots + v_{i-1}}(m) > f^{v_1+ \cdots  +v_{i-1}+1}(m) > \cdots > f^{v_1+ \cdots  +v_{i-1}+v_i}(m).
\eeq
The arithmetic function $f$ is \emph{wildly increasing-decreasing} if,  
for every finite sequence $V$ of positive integers, there exists an integer $m \in \Omega$ 
such that $m$ has increasing-decreasing pattern $V$ with respect to $f$.

More generally, let $\mcf = (f_j)_{j=0}^{\infty}$ be a sequence of arithmetic functions 
$f_j:\Omega \rightarrow \Omega$.  The \mcf-trajectory of $m \in \Omega$ is 
the sequence $\left(f_j(m) \right)_{j=1}^{\infty}$. 
Let $V = (v_i)_{i=1}^n$ be a finite sequence of positive integers.  
We say that an integer $m$ in $\Omega$ has \emph{increasing-decreasing pattern}\index{pattern} 
$V$ with respect to the sequence $\mcf$ if, for odd $i$, 
\[
f_{v_1+ \cdots + v_{i-1}}(m) < f_{v_1+ \cdots + v_{i-1}+1}(m) < \cdots < f_{v_1+ \cdots + v_{i-1}+v_{i}}(m)
\]
and, for even $i$, 
\[
f_{v_1+ \cdots + v_{i-1}}(m) > f_{v_1+ \cdots  +v_{i-1}+1}(m) > \cdots > f_{v_1+ \cdots  +v_{i-1}+v_i}(m).
\]
The sequence \mcf\ is \emph{wildly increasing-decreasing} if,  
for every finite sequence $V$ of positive integers, there exists $m \in \Omega$ 
such that $m$ has increasing-decreasing pattern $V$ with respect to the sequence \mcf.  

Fix an arithmetic function $f:\Omega \rightarrow \Omega$.
Let $(k_j)_{j=0}^{\infty}$ be a strictly increasing sequence of nonnegative integers.   
If $f$ is wildly increasing-decreasing, then the sequence of arithmetic functions 
$\mcf = \left( f^{k_j} \right)_{j=0}^{\infty}$ is also wildly increasing-decreasing.

It is not known if aliquot sequences are wildly increasing-decreasing.  
Lenstra~\cite{lens75} proved that,  for every positive integer $v_1$,   
there are infinite many integers $m$ 
such that $m < s(m) < s^2(m) < \cdots < s^{v_1}(m)$.   
Erd\H os~\cite{erdo76} subsequently refined this result.  
Pomerance~\cite{pome20} proved that,  for every positive integer $v_2$, 
there are infinitely many integers $m$ 
such that $m > s(m) > s^2(m) > \cdots > s^{v_2}(m)$. 
It is an open problem to determine if, for every pair of positive integers $v_1, v_2$ 
there are integers $m$ such that 
\[
m < s(m) < s^2(m) < \cdots < s^{v_1}(m) 
\]
and 
\[
s^{v_1}(m) > s^{v_1+1}(m) > s^{v_1+2}(m) > \cdots > s^{v_1+v_2}(m).
\]

 There are two 
competing conjectures about the trajectories of the sum of the proper divisors function $s(m)$.  
Catalan~\cite{cat88} and Dickson~\cite{dic13}  
conjectured that every aliquot sequence is bounded.  
Guy and Selfridge~\cite{guy-sel71,guy-sel75}  conjectured that infinitely 
many aliquot sequences go to infinity.  

The goal of this paper is to prove that Syracuse function $S(m)$ 
defined by~\eqref{Syracuse:function} is wildly increasing-decreasing.   
 The proof uses integer matrices and systems of linear diophantine equations.

\section{Integer matrices}        \label{Syracuse:section:IntegerMatrices}

We consider  vectors in \Rn.
A  vector is \emph{positive}\index{positive vector}\index{vector!positive}  
if all of its coordinates are positive and  \emph{negative}\index{negative vector}\index{vector!negative}  
if all of its coordinates are negative.
A vector is \emph{integral}\index{integral vector}\index{vector!integral}  
if all of its coordinates are integers.  

Let $\Z^n$ denote the set of $n$-dimensional integral vectors. 
An integral vector is \emph{odd}\index{odd vector}\index{vector!odd} 
if its coordinates are odd integers and \emph{even}\index{even vector}\index{vector!even} 
if its coordinates are even integers. 
An integral vector is \emph{primitive}\index{primitive vector}\index{vector!primitive} 
if its coordinates are relatively prime (not necessarily pairwise relatively prime) positive integers.

\bt                         \label{Syracuse:theorem:M}
Let $(a_i)_{i=1}^n$ and $(b_j)_{j=2}^{n+1}$ be sequences of nonzero real numbers.  
Consider the $n \times (n+1)$ matrix 
\[
M = \bmat 
a_1 & -b_2 & 0          & 0          & 0             & 0        &  \cdots          & 0 & 0 & 0 \\
0         & a_2  & -b_3  & 0          & 0             & 0      &  \cdots           & 0& 0 & 0  \\
0         & 0           & a_3  & -b_4 & 0             & 0      &  \cdots           & 0 & 0 & 0 \\
0         & 0           & 0            & a_4 & -b_5  & 0       &  \cdots           & 0 & 0 & 0 \\
\vdots & &&&&&&&& \vdots \\
0         & 0           & 0            &0          & 0             & 0       & \cdots            & 0 & a_n & -b_{n+1} \\
\emat
\]
\benum
\item[(a)]
The matrix $M$ has rank $n$ and the kernel of $M$ has dimension 1. 

\item[(b)]
Let $(a_i)_{i=1}^n$ and $(b_j)_{j=2}^{n+1}$ be sequences of positive real numbers.  
If $\mbx  \in \kernel(M)$ and $\mbx \neq \mbo$, 
then the vector \mbx\ is positive or negative.   

\item[(c)]
Let $(a_i)_{i=1}^n$ and $(b_j)_{j=2}^{n+1}$ be sequences of positive integers. 
There is a unique primitive vector $\mbz$ that generates $\kernel(M)$.  

\item[(d)]
Let $(a_i)_{i=1}^n$ be a sequence of odd positive integers  
and let  $(b_j)_{j=2}^{n+1}$ be a sequence of even positive integers. 
If $\mbz = \vectorsmallzn+1$ is the unique 
primitive vector in $\kernel(M)$,  then the integers $z_i$ are even for all $i \in \{1,2,\ldots, n\}$ 
and the integer $z_{n+1}$ is odd. 

\item[(e)]
Let $(a_i)_{i=1}^n$ be a sequence of odd positive integers  
and let  $(b_j)_{j=2}^{n+1}$ be a sequence of even positive integers. 
Let $\mbh \in \Z^n$ be an odd vector.  
If $M\mbx = \mbh$ for some integral vector $\mbx \in \Z^{n+1}$, 
then there exists an odd positive vector $\mbg \in \Z^{n+1}$ such that  
\[
M\mbg = \mbh.
\]
\eenum
\et

\begin{proof} \

\benum
\item[(a)]
The $n$ row vectors of the matrix $M$ are linearly independent  and so $M$ has rank $n$. 

\item[(b)] 
Let $\mbx = \vectorsmallxn+1 \in \kernel(M)$.  
If $M \mbx = \mbo$, then for all $i \in \{1,\ldots, n\}$ we have 
\[
a_i x_i = b_{i+1} x_{i+1}.
\]
If $\mbx \neq \mbo$, then $x_{i_0} \neq 0$ for some ${i_0} \in \{1,\ldots, n+1\}$.  
The positivity of the numbers $a_i$ and $b_i$ implies that if $x_{i_0} > 0$, 
then $x_i > 0$ for all $i \in \{1,\ldots, n+1\}$.  
Similarly, if $x_{i_0} < 0$, then $x_i < 0$ for all $i \in \{1,\ldots, n+1\}$.  
Therefore, \mbx\ is either a positive vector or a negative vector.

\item[(c)]
Consider $M$ as a matrix with coordinates in the field \Q\ of rational numbers.  
The kernel of $M$ is one-dimensional.  Let $\mbx \in \Q^{n+1}$ be a nonzero vector 
in the kernel of $M$.  By~(b), the coordinates of \mbx\ are either all positive rational 
numbers or all negative rational numbers.  
Multiplying, if necessary, by -1, we can assume that the coordinates are positive rational numbers.
Multiplying by a common denominator of the denominators of the coordinates, 
we obtain a vector whose coordinates are positive integers.  Dividing by the greatest common 
divisor of these integers, we obtain a primitive vector in the kernel of $M$.  

Let $\mbz  = \vectorsmallzn+1$ 
and $\mbz' = \vectorsmallz'n+1$
be primitive vectors in $\kernel(M)$. 
Because the kernel is one-dimensional, there is a positive rational number $p/q$ 
with $\gcd(p,q) = 1$ such that $(p/q)\mbz = \mbz'$.  Multiplying by $q$ we obtain 
\[
\bmat pz_1 \\ \vdots \\ pz_n \\ pz_{n+1} \emat = p\mbz 
= q \mbz' = \bmat qz'_1 \\ \vdots \\ qz'_n \\ qz'_{n+1} \emat 
\]  
and so 
\[
pz_i = qz'_i 
\]
for all $i \in \{1,\ldots, n, n+1\}$. 
The divisibility condition $\gcd(p,q) = 1$ implies that $p$ divides $z'_i$ and 
$q$ divides $z_i$ for all $i \in \{1,\ldots, n, n+1\}$.  Because the vectors \mbz\ and \mbz' are primitive, 
we have $p=q=1$ and so $\mbz = \mbz'$.  Thus, $\kernel(M)$ contains a unique primitive vector.

\item[(d)]
Let $\mbz  = \vectorsmallzn+1$ be the  unique primitive vector in $\kernel(M)$. 
For all $i \in \{1,\ldots, n\}$ we have 
\[
a_iz_i = b_{i+1} z_{i+1}.
\]
Because $b_i$ is even and $a_i$ is odd, it follows that 
$z_i$ is even for all $i \in \{1,\ldots, n\}$.  
We have $\gcd(z_1,\ldots, z_n, z_{n+1}) = 1$ 
because the vector $\mbz$ is primitive, 
and so $z_{n+1}$ must be odd. 

\item[(e)] 
Let $\mbh = \vectorsmallhn \in \Z^n$ be an odd vector.    
If $\mbx = \vectorsmallxn+1 \in \Z^{n+1}$ and $M\mbx = \mbh$, then 
for all $i \in \{1,\ldots, n\}$ we have 
\[
a_i x_i - b_{i+1} x_{i+1} = h_i. 
\] 
Because the integer $b_{i+1}$ is even and the integers $a_i$ and $h_i$ are odd, it follows that 
the integer $x_i$ is odd for all $i \in \{1,\ldots, n\}$.  
The integer $x_{n+1}$ is not necessarily odd, nor are the integers 
$x_1,\ldots, x_n, x_{n+1}$ necessarily positive.   

Let $\mbz = \vectorsmallzn+1$ be the unique 
primitive vector in $\kernel(M)$. 
By~(d), the integers $z_i$ are even for all $i \in \{1,2,\ldots, n\}$ and 
the integer $z_{n+1}$ is odd. 
We have $M\mbz = \mbo$ and so 
\[
M(\mbx + k\mbz) = M\mbx = \mbh
\]
for all integers $k$, where 
\[
\mbx + k\mbz = \vectorsmallxn+1 + k\vectorsmallzn+1 
= \bsmallmat x_1 + kz_1 \\ \vdots \\ x_n + kz_n \\ x_{n+1}+kz_{n+1} \esmallmat.
\]
The coordinates $z_i$ are positive, and so the vector $\mbx+ k\mbz$ is positive 
for all sufficiently large $k$.  

For all $i \in \{1,\ldots, n\}$ the integer $x_i$ is odd and the integer $z_i$ is even, 
and so $x_i+kz_i$ is an odd integer.  
The integer $z_{n+1}$ is odd, and so the integers $x_{n+1} +kz_{n+1}$ 
and $x_{n+1}+(k+1)z_{n+1}$ have opposite parity.  
It follows that for all sufficiently large $k$, either $\mbg = \mbx+k\mbz$ 
or \mbg = $\mbx + (k+1)\mbz$ is an odd positive vector such that $M\mbg = \mbh$.  
\eenum
This completes the proof.  
\end{proof}

\bt                 \label{Syracuse:theorem:Mab}
Let $(a_i)_{i=1}^n$ and  $(b_j)_{j=2}^{n+1}$  
be sequences of nonzero  integers such that 
\beq            \label{Syracuse:divisibility}
\gcd(a_i,b_j) = 1
\eeq
for all $i \in \{1,\ldots, n\}$ and $j \in \{2,\ldots, n+1\}$.  
Let $M$ be the $n \times (n+1)$ matrix defined 
in Theorem~\ref{Syracuse:theorem:M}.    
The homomorphism $M: \Z^{n+1} \rightarrow \Z^n$ is surjective. 
\et

\begin{proof}
We give two proofs.  
The first proof uses induction on $n$.  
Because $(a_1,b_2)=1$, for every integer $h_1$ there are integers $x_1,x_2$ such 
that $a_1x_1 + b_2x_2 = h_1$.  
Equivalently, for the $2 \times 1$ matrix $M = \bmat a_1 & b_2 \emat$ 
and for $\mbh = (h_1)$, we have 
\[
M \bmat x_1 \\ x_2 \emat  = \bmat a_1 & b_2 \emat \bmat x_1 \\ x_2 \emat =   \mbh.
\]
This is the case $n=1$. 

For $n \geq 2$, let  $M$ be the $n \times (n+1)$ matrix defined 
in Theorem~\ref{Syracuse:theorem:M} and let $\mbh = \vectorsmallhn \in \Z^n$. 
There exists an integral vector $\mbx = \vectorsmallxn+1 \in \Z^{n+1}$ 
such that 
\[
M\mbx = \mbh 
\]
if and only if the diophantine system of $n$ linear equations in $n+1$ variables 
\begin{align*}
a_1 x_1 - b_2 x_2 & = h_1 \\ 
a_2 x_2 - b_3 x_3 & = h_2 \\ 
\vdots & \\
a_{i-1} x_{i-1} - b_{i}x_i & = h_{i-1}\\ 
a_i x_i - b_{i+1}x_{i+1} & = h_i \\ 
\vdots \\
a_n x_n - b_{n+1}x_{n+1} & = h_n 
\end{align*}
has a solution in integers $x_1,\ldots, x_n, x_{n+1}$.

For $n=2$, we have two equations  
\begin{align*}
a_1 x_1 - b_2x_2 & = h_1 \\ 
a_2 x_2 - b_3x_3 & = h_2.
\end{align*}
The divisibility condition~\eqref{Syracuse:divisibility} 
gives $\gcd(a_1,b_2) = \gcd(a_2,b_3) = 1$, and so both equations have 
 solutions in integers.  If $(c_1,d_2)$ is a particular solution  
of the first equation, then the general solution of the first equation is 
\begin{align*}
x_1 & = c_1 + b_2y_1 \\
x_2 & = d_2 + a_1y_1
\end{align*}
for any integer $y_1$.  
 If $(c_2,d_3)$ is a particular solution  
of the second equation, then the general solution of the second equation is 
\begin{align*}
x_2 & = c_2 + b_3y_2 \\
x_3 & = d_3 + a_2y_2
\end{align*}
for any integer $y_2$.  
We have a simultaneous solution of the system of two equations if and only if there exist integers 
$y_1$ and $y_2$ such that
\[
x_2 = d_2+a_1y_1 = c_2+b_3y_2
\]
or, equivalently, 
\[
a_1y_1 - b_3y_2 = c_2-d_2.
\]
Because $\gcd(a_1,b_3) = 1$, this equation has a solution in integers.  
This proves the Theorem for $n=2$.

Let $n \geq 3$ and assume that the Theorem is true for $n-1$ equations 
in $n$ variables.    
Consider the diophantine system of $n$ equations in $n+1$ variables 
such that, for all $i \in \{1,\ldots, n\}$,    equation  $(i)$ is 
\[
a_i x_i - b_{i+1}x_{i+1}  = h_i. 
\]
This equation has a solution in integers because $\gcd(a_i,b_{i+1}) = 1$.  
If $(c_i,d_{i+1})$ is a particular solution of the equation, 
then the general solution of equation  $(i)$ is 
\begin{align*}
x_i & = c_i + b_{i+1} y_i \\
x_{i+1} & = d_{i+1} + a_i y_i
\end{align*}
for any integer $y_i$. 
Similarly, for all $i \in \{2,\ldots, n+1\}$,   equation   $(i-1)$ is 
\[
a_{i-1} x_{i-1} - b_ix_{i}  = h_{i-1} 
\]
This equation has a solution in integers because $\gcd(a_{i-1},b_i) = 1$.  
If $(c_{i-1},d_i)$ is a particular solution of the equation, 
then the general solution of  equation   $(i-1)$  is 
\begin{align*}
x_{i-1} & = c_{i-1} + b_iy_{i-1} \\
x_i & = d_i + a_{i-1} y_{i-1}
\end{align*}
for any integer $y_{i-1}$. Equations $(i-1)$ and   $(i)$ have a simultaneous solution 
in integers if and only if there exist integers $y_{i-1}$ and $y_i$ such that 
\[
d_i + a_{i-1} y_{i-1} = c_i + b_{i+1} y_i 
\]
or, equivalently, if 
\[
a_{i-1} y_{i-1} -  b_{i+1}y_i = c_i - d_{i-1}. 
\]
It follows that the original system of $n$ equations in $n+1$ variables 
has a solution in integers if and only if the following  system of $n-1$ equations 
in $n$ variables  has a solution in integers: 
\begin{align*}
a_{1} y_1 -  b_3y_2& = c_2 - d_1 \\
a_2 y_2-  b_4y_3 & = c_3 - d_2 \\
\vdots & \\
a_{n-1} y_{n-1} -  b_{n+1}y_n & = c_n - d_{n-1}. 
\end{align*}
The divisibility condition~\eqref{Syracuse:divisibility} and the induction hypothesis 
imply that this system 
of equations does have an integral solution.  This completes the first proof. 

The second proof uses the \emph{Smith normal form} of an integral matrix 
(Marcus and Minc~\cite[pp. 40--48]{marc-minc92}).  
The Smith normal form of an $m \times n$ matrix $M$ of rank $k$ is the unique diagonal matrix 
\[
SNF(M) = \bmat 
s_1 & 0 & 0          & 0          & 0                     &  \cdots          & 0 & 0 & 0 \\
0         & s_2  & 0  & 0          & 0                &  \cdots           & 0& 0 & 0  \\
0         & 0           & s_3 & 0 & 0             &  \cdots           & 0 & 0 & 0 \\
\vdots & &&&&&&&  \vdots \\ 
0         & 0           & 0            & 0 & 0      &  \cdots           & s_k & 0 & 0 \\
0         & 0           & 0            &0          & 0          & \cdots            & 0 & 0 &0 \\
\vdots & &&&&&&& \vdots \\
0         & 0           & 0            &0          & 0            & \cdots            & 0 & 0 &0 \\
\emat 
\]
constructed from $M$ by integral elementary row and column operations 
and whose coordinates $s_1,\ldots, s_k$ are positive integers. 
The \emph{$k$th determinantal divisor} of $A$ is the greatest common divisor 
of all of the $k \times k$ minors of $A$.  
If the $k$th determinantal divisor of $A$ is $1$, then $s_i = 1$ for all $i = 1,\ldots, k$.
If $A$ has rank $n$ and the $n$th determinantal divisor of $A$ is $1$, 
then $s_i = 1$ for all $i = 1,\ldots, n$ and 
\beq                          \label{Syracuse:SNF}
SNF(M) = \bmat 
1 & 0 & 0          & 0          & 0        &  \cdots          & 0 & 0 & 0 \\
0         & 1  & 0  & 0          & 0          &  \cdots           & 0& 0 & 0  \\
0         & 0           & 1  & 0 & 0           &  \cdots           & 0 & 0 & 0 \\
\vdots & &&&&&&& \vdots \\
0         & 0           & 0            &0           & 0       & \cdots            & 0 & 1 &0 \\
\emat. 
\eeq
The homomorphism $M: \Z^{n+1} \rightarrow \Z^n$ 
is surjective if and only if the  homomorphism $SNF(M): \Z^{n+1} \rightarrow \Z^n$ 
is surjective if and only if  $\rank(M) = n$ and $s_i = 1$ for all $i = 1,\ldots, n$.

Let $M$ be the $n \times (n+1)$ matrix constructed  
in Theorem~\ref{Syracuse:theorem:M}  from integer sequences 
 $(a_i)_{i=1}^n$ and  $(b_j)_{j=2}^{n+1}$ such that 
$\gcd(a_i,b_j) = 1$
for all $i \in \{1,\ldots, n\}$ and $j \in \{2,\ldots, n+1\}$.
The determinant of the $n \times n$ minor obtained by deleting the first column of $M$
is $(-1)^n \prod_{j=2}^{n+1} b_j$.  
The determinant of the $n \times n$ minor obtained by deleting the last  column of $M$  
is $ \prod_{i=1}^{n} a_i$. 
The divisibility condition $\gcd(a_i,b_j) = 1$ implies that  
these determinants are relatively prime integers and so the $n$th determinantal divisor of $M$ is 1. 
 It follows that the matrix $M$ has Smith normal form~\eqref{Syracuse:SNF}  
and the homomorphisms $SNF(M):\Z^{n+1} \rightarrow \Z^n$ and $M:\Z^{n+1} \rightarrow \Z^n$ 
are surjective. 
This completes the second proof. 
\end{proof}

\bt                               \label{Syracuse:theorem:M34}
Let $(v_i)_{i=1}^{n+1}$ be a sequence of positive integers.  
Define the positive integral sequences 
$(a_i)_{i=1}^n$ and $(b_j)_{j=2}^{n+1}$ as follows:
\[
a_i = 3^{v_i}
\]
and 
\[
b_j = \begin{cases}
4^{v_j} & \text{if $j$ is even} \\ 
2^{v_j} & \text{if $j$ is odd.}
\end{cases}
\]
Let 
\[
c_{n+1} = \begin{cases}
4 & \text{if $n+1$ is even} \\ 
2 & \text{if $n+1$ is odd.}
\end{cases}
\]
The system of linear diophantine equations
\begin{align*}
3^{v_1} x_1 - 4^{v_2} x_2 & = 1 \\
3^{v_2} x_2 - 2^{v_3} x_3 & = - 1 \\
3^{v_3} x_3 - 4^{v_4} x_4 & =  1 \\
\vdots & \\
3^{v_n} x_n - c_{n+1}^{v_{n+1}} x_{n+1} & = (-1)^{n+1} 
\end{align*} 
has a solution in odd positive integers $w_1, w_2, \ldots, w_n, w_{n+1}$. 
\et

\begin{proof}
Let 
\[
M = \bmat 
3^{v_1} & -4^{v_2} & 0          & 0          & 0             & 0        &  \cdots          & 0 & 0 & 0 \\
0         & 3^{v_2}  & -2^{v_3}  & 0          & 0             & 0      &  \cdots           & 0& 0 & 0  \\
0         & 0           & 3^{v_3}  & -4^{v_4} & 0             & 0      &  \cdots           & 0 & 0 & 0 \\
0         & 0           & 0            & 3^{v_4} & -2^{v_5}  & 0       &  \cdots           & 0 & 0 & 0 \\
\vdots & &&&&&&&& \vdots \\
0         & 0           & 0            &0          & 0             & 0       & \cdots            & 0 & 3^{v_n} & -c_{n+1}^{v_{n+1}} 
\emat
\]
 be the $n \times (n+1)$ matrix constructed  
in Theorem~\ref{Syracuse:theorem:M}. 
Let $h_i = (-1)^{i+1}$ for $i \in \{1,\ldots, n\}$ and  $\mbh = \vectorsmallhn \in \Z^n$.  
By  Theorem~\ref{Syracuse:theorem:Mab}, 
there exists an integralr vector $\mbx\in \Z^{n+1}$ such that $M \mbx = \mbh$. 
By Theorem~\ref{Syracuse:theorem:M}(e), 
there exists an odd positive vector $\mbw = \vectorsmallwn+1\in \Z^{n+1}$ such that 
$M \mbw = \mbh$.
This completes the proof. 
\end{proof}

\section{Iterations of the Syracuse function}              \label{Syracuse:section:proofs}

Let $\Omega$ be  the set of odd positive integers. 
The Syracuse function  is the arithmetic function $S:\Omega \rightarrow \Omega$  
defined by 
\[
S(m) = \frac{3m + 1}{2^e}
\]
where $e$ is the largest integer such that $2^e$ divides  $3m+1$. 
Note that $S(m) = 1$ if and only if $m = (4^e-1)/3$ for some positive integer $e$.

\bl                                                                   \label{Syracuse:lemma:SJM-up}  
Let $m$ be a positive integer such that 
\[
m \equiv -1 \pmod{4} 
\]
and let $v$ and $w$ be the unique positive integers with $w$ odd such that 
\[
m = 2^{v+1} w  - 1.
\] 
Then 
\beq              \label{Syracuse:SJM-up-2}
m < S(m) < S^2(m) < \cdots < S^{v}(m) 
\eeq
and 
\beq              \label{Syracuse:SJM-up-end}
 S^v(m) =  2\cdot 3^{v}  w - 1.
\eeq
Moreover, 
\[
S^{v+1}(m) = \frac{3^{v} w -1}{2^e}
\]
for some positive integer $e$. 
\el

\begin{proof} 

We use induction on $j$ to prove that if 
\[
0 \leq j \leq v
\]
then 
\beq              \label{Syracuse:SJM-up-1}
S^j(m) = 2^{v+1-j}  3^j  w - 1 = \left( \frac{3}{2}\right)^j (m+1) - 1.
\eeq
This is true for $j=0$.  For $j=1$ we have 
\[
S(m) = \frac{3 (2^{v+1} w-1)+1}{2} = \left( \frac{3}{2} \right) 2^{v+1}w -1 =   \left( \frac{3}{2} \right) (m+1) -1. 
\]
If $j \leq v-1$ and 
\[
S^j(m) = 2^{v+1-j}  3^j  w - 1
\]
then 
\[
S^{j+1}(m) = 2^{v+1-(j+1)}3^{j+1}w-1 = \left( \frac{3}{2}\right)^{j+1} (m+1) - 1. 
\]
This proves~\eqref{Syracuse:SJM-up-1}, and~\eqref{Syracuse:SJM-up-1} 
implies~\eqref{Syracuse:SJM-up-2}. 

Note that  
\[
S^{v}(m) = 2\cdot 3^{v}  w - 1
\]
implies 
\[
S^{v+1}(m) = S\left( 2 \cdot 3^{v} w - 1 \right) = \frac{3^{v+1}w -1}{2^e} 
\] 
for some positive integer $e$. 
This completes the proof.  
\end{proof}

\bl                                                                     \label{Syracuse:lemma:SJM-down}
Let   $m$ be a positive integer such that 
\[
m \equiv 1 \pmod{8} 
\]
and let $v$ and $w$ be the unique positive integers  with $v \geq 3$ and  $w$ odd such that  
\[
m = 2^v w+1.
\]
If 
\[
v_0 = \left[ \frac{v-1}{2} \right]  
\]
then 
\[
v = 2v_0+1 \quad\text{ or } \quad v = 2v_0+2
\]
and 
\beq              \label{Syracuse:SJM-down-2}
m > S(m) > S^2(m) > \cdots > S^{v_0-1}(m) > S^{v_0}(m) .
\eeq
If $v = 2v_0+1$, then 
\beq              \label{Syracuse:SJM-down-2-end} 
S^{v_0}(m)  = 2  \cdot 3^{v_0} w + 1
\eeq 
If $v = 2v_0+2$, then 
 \[
S^{v_0}(m)  = 4 \cdot 3^{v_0} w + 1
\]
\el

\begin{proof}
We have $v-2j \geq 1$  if and only if 
\[
j \leq v_0 = \left[\frac{v-1}{2}\right]
\] 
and  $v \geq 3$ implies $v_0 \geq 1$. 
We use induction on $j$ to prove that if 
\[
0 \leq j \leq v_0 
\]
then 
\beq              \label{Syracuse:SJM-down-1}
S^j(m) = 2^{v-2j} 3^j w + 1 =  \left( \frac{3}{4}\right)^j (m-1)+1. 
\eeq
This is true for $j = 0$.  For $j=1$ we have 
\begin{align*}
S(m) & = S\left( 2^vw+1\right) = \frac{3(2^vw+1)+1}{4}  \\ 
& = 2^{v-2}3w +1 = \frac{3}{4} (2^v w)+1 \\ 
&  = \frac{3}{4} (m-1)+1.  
\end{align*}

If $1 \leq j \leq v_0-1$ and 
\[
S^j(m) = 2^{v-2j} 3^j w + 1
\]
then 
\[
v - 2(j+1) \geq v - 2v_0 = v - 2\left[\frac{v-1}{2}\right] \geq 1 
\]
and 
\begin{align*}
S^{j+1}(m) & = S\left( 2^{v-2j} 3^j w + 1\right)  = \frac{  2^{v-2j} 3^{j+1} w + 4}{4} \\
& =  2^{v-2(j+1)} 3^{j+1} w +1 = \left(\frac{3}{4}\right)^{j+1} (2^v w)+1 \\
& =  \left(\frac{3}{4}\right)^{j+1} (m-1)+1.  
\end{align*}
This proves~\eqref{Syracuse:SJM-down-1}, and~\eqref{Syracuse:SJM-down-1} 
implies~\eqref{Syracuse:SJM-down-2}. 

We have 
\[
S^{v_0}(m)  = 2^{v-2v_0} 3^{v_0} w + 1
\]
If $v = 2v_0+1$, then 
 \[
S^{v_0}(m)  = 2  \cdot 3^{v_0} w + 1
\] 
If $v = 2v_0+2$, then 
 \[
S^{v_0}(m)  = 4 \cdot 3^{v_0} w + 1
\]
This completes the proof. 
\end{proof}

\bt              \label{Syracuse:theorem:wildly}
The Syracuse function $S$ is wildly increasing-decreasing. 
\et

\begin{proof} 
Let 
\[
c_{n+1} = \begin{cases}
4 & \text{if $n+1$ is even} \\ 
2 & \text{if $n+1$ is odd.}
\end{cases}
\]

Let $V = (v_i)_{i=1}^n$ be a finite sequence of positive integers.  
By Lemma~\ref{Syracuse:lemma:SJM-up}, if 
\[
m = 2^{v_1+1} w_1 - 1
\]
for some odd positive integer $w_1$, then the odd positive integer $m$ 
increases for $v_1$ iterations of the Syracuse function $S$.  
From~\eqref{Syracuse:SJM-up-end} we have  the odd positive integer  
\[
S^{v_1}(m) = 2\cdot 3^{v_1} w_1 - 1.
\]

By Lemma~\ref{Syracuse:lemma:SJM-down}, 
if 
\[
S^{v_1}(m)  = 2\cdot 4^{v_2} w_2+1
\]
for some  odd positive  integer $w_2$, then 
the odd positive integer $S^{v_1}(m)$ decreases for $v_2$ iterations 
of the Syracuse function $S$.   
The integers $w_1$ and $w_2$ are solutions of the diophantine equation 
\[
2\cdot 3^{v_1} w_1 - 1 = 2\cdot 4^{v_2} w_2+1
\]
or, equivalently, 
\[
3^{v_1} w_1 - 4^{v_2} w_2 = 1.
\]
From~\eqref{Syracuse:SJM-down-2-end} we have  the odd positive integer 
\[
S^{v_1+v_2}(m) = 2\cdot 3^{v_2} w_2 + 1.
\]

By Lemma~\ref{Syracuse:lemma:SJM-up}, if 
\[
S^{v_1+v_2}(m) =   2^{v_3 +1} w_3 - 1
\]
for some odd positive integer $w_3$, 
then the integer $S^{v_1+v_2}(m)$ increases for $v_3$ iterations 
of the Syracuse function $S$.  
The integers $w_2$ and $w_3$ are solutions of the diophantine equation 
\[
2\cdot 3^{v_2} w_2 + 1 =  2^{v_3 +1} w_3 - 1
\]
or, equivalently, 
\[
3^{v_2} w_2 - 2^{v_3} w_3 = - 1.
\]
Continuing inductively, we obtain an odd positive  integer $m$ 
that has increasing-decreasing pattern $V$ with respect to the 
Syracuse function $S$ if the system of linear diophantine equations 
\begin{align*}
3^{v_1} w_1 - 4^{v_2} w_2 & = 1 \\
3^{v_2} w_2 - 2^{v_3} w_3 & = - 1 \\
3^{v_3} w_3 - 4^{v_4} w_4 & = 1 \\
\vdots & \\
3^{v_n} w_n - c_{n+1}^{v_{n+1}} w_{n+1} & = (-1)^{n+1} 
\end{align*} 
has a solution $w_1, \ldots, w_{n+1}$ in odd positive integers.   
This is precisely what Theorem~\ref{Syracuse:theorem:Mab} provides.
This completes the proof. 
\end{proof}

\section{A simpler proof}
I thank an anonymous referee for the reference to Everett~\cite{ever77} 
and the following variant of Theorem~\ref{Syracuse:theorem:wildly}.  

\bt            \label{Syracuse:theorem:simpler}
Let $k$ be a positive integer and let $\Lambda = (\lambda_1,\lambda_2,\ldots, \lambda_k)$ be a $k$-tuple 
with $\lambda_i \in \{\text{increase}, \text{decrease}\}$ for all $i \in \{1,2,\ldots, k\}$.  
Let $S(\Lambda)$ be the set of all odd positive integers $m$ such that 
\[
S^{i-1}(m)  < S^i(m) \qquad \text{if $\lambda_i = $} increase  
\] 
and 
\[
S^{i-1}(m) > S^i(m) \qquad \text{if $\lambda_i = $} decrease    
\] 
 for all $i \in \{1,2,\ldots, k\}$.  
There are positive integers  $\ell = \ell(\Lambda)$ and 
$r = r(\Lambda)$ with $r$ odd 
such that, for all $m \geq 1$, 
\beq            \label{Syracuse:simple}
m \equiv r \pmod{2^{\ell}} \qquad \text{implies} \qquad  m \in S(\Lambda).
\eeq
\et

This result immediately implies thats that the Sylvester function is wildly increasing-decreasing.

\begin{proof}
The proof is by induction on $k$. 

Let $k=1$.  
If $\Lambda = \{\text{increase}\}$, then choose $\ell = 2$ and $r=3$.  
By Lemma~\ref{Syracuse:lemma:SJM-up}, the congruence 
$m \equiv 3 \pmod{2^2}$ implies $m \in S(\Lambda)$. 
If $\Lambda = \{\text{decrease}\}$, then choose $\ell = 3$ and $r=1$.  
By Lemma~\ref{Syracuse:lemma:SJM-down}, the congruence  
$m \equiv 1 \pmod{2^3}$  implies $m \in S(\Lambda)$.
This proves the Theorem for $k=1$.

Let $k \geq 1$ and assume that for every $k$-tuple $\Lambda$ there exist 
positive integers  $\ell = \ell(\Lambda)$ and 
$r = r(\Lambda)$ with $r$  odd 
that satisfy~\eqref{Syracuse:simple} for all $m \geq 1$. 
Let $\Lambda' = (\lambda_1,\lambda_2,\ldots, \lambda_k, \lambda_{k+1})$ be a $(k+1)$-tuple 
with $\lambda_i \in \{\text{increase}, \text{decrease}\}$ for all $i \in \{1,2,\ldots, k, k+1\}$.  
Consider the $k$-tuple $\Lambda = (\lambda_2,\ldots, \lambda_k, \lambda_{k+1})$.  
By the induction hypothesis, there exist positive integers $\ell'$ and $r'$ 
with $r'$ odd such that, for all $m' \geq 1$, 
\[
m' \equiv r' \pmod{2^{\ell'}} \qquad \text{implies} \qquad  m' \in S(\Lambda').
\]

Suppose $\lambda_1 = \text{increase}$.  The odd positive integer $m$ is in $S(\Lambda)$ 
if 
\[
S(m) = \frac{3m+1}{2}  \equiv r' \pmod{2^{\ell'}}  
\]
or, equivalently, if 
\[
3m \equiv 2r'-1 \pmod{2^{\ell' + 1}}.
\]
Let $\ell = \ell' + 1$.  
Choose $t \in \N$ such that $3t \equiv 1 \pmod{2^{\ell}}$.  Then $t$ and $r = (2r'-1)t$ are odd and $\ell$ and  $r$ satisfy~\eqref{Syracuse:simple}.  
 
 Suppose $\lambda_1 = \text{decrease}$.  The odd positive integer $m$ is in $S(\Lambda)$ 
if 
\[
S(m) = \frac{3m+1}{4}  \equiv r' \pmod{2^{\ell'}}   
\]
or, equivalently, if 
\[
3m \equiv 4r'-1 \pmod{2^{\ell' + 2}}.
\]
Let $\ell = \ell' + 2$.  
Choose $t \in \N$ such that $3t \equiv 1 \pmod{2^{\ell}}$.  
Then $t$ and $r = (4r'-1)t$ are odd and  
$\ell$ and  $r$ satisfy~\eqref{Syracuse:simple}.  
Thus, the Theorem holds for all $(k+1)$-tuples $\Lambda$. 
This completes the proof. 
\end{proof}

\section{Open problems} 
\subsection{The infinite constant sequence $v_i=1$}
The Collatz conjecture states that for every positive integer $m$ 
there is an integer $k_m$ such that $S^{k_m}(m) = 1$ 
and so the Syracuse trajectory 
$\left( S^j(m)\right)_{j=0}^{\infty}$ is eventually constant.  
This would imply that there is no infinite sequence of positive integers  $V = \left( v_i \right)_{i=1}^{\infty}$ 
and no positive integer $m$ for which the trajectory $\left( S^j(m) \right)_{j=0}^{\infty}$ 
satisfies the increasing-decreasing 
conditions~\eqref{Syracuse:increasing} and ~\eqref{Syracuse:decreasing} for all $i \in \N$. 
Is it possible to prove that some particular infinite sequence 
$V = \left( v_i \right)_{i=1}^{\infty}$ is not 
the increasing-decreasing pattern of any positive integer $m$ 
under iterations of the Syracuse function?  For example, 
can the constant sequence  $V = \left( v_i \right)_{i=1}^{\infty}$ with $v_i = 1$ for all $i$ 
be proven impossible? This is equivalent to proving that if $M$ 
is the infinite matrix with $(i,j)$th coordinate 
\[
M_{i,j} = \begin{cases}
3 & \text{if $j = i$}\\
-4 & \text{if $j = i+1$} \\
0 & \text{otherwise} 
\end{cases}
\]
for all $i,j \in \N$ and if $\mbh =  (h_i)_{i=1}^{\infty}$ is the infinite vector 
with $h_i = (-1)^{i+1}$ for all $i \in \N$, 
then there exists no infinite odd positive vector $\mbw =  (w_j)_{j=1}^{\infty}$ such that $M\mbw = \mbh$. 

\subsection{The finite constant sequence $v_i=1$}
For every positive integer $n$ there exists an odd positive integer $m$ 
that has the increasing-decreasing pattern
 $V = \left( v_i \right)_{i=1}^n$ with $v_i = 1$ for all $i \in \{1,2,\ldots, n\}$. 
 Let $m_n$ be the least such integer.  
The sequence $(m_n)_{n=1}^{\infty}$ is increasing.  Is it strictly increasing? 
Can $m_n$ be efficiently computed?

\subsection{Sufficient conditions for wildly increasing-decreasing} 
The method of proof of Theorem~\ref{Syracuse:theorem:simpler} 
can be extended to prove that a large class of arithmetic functions 
is wildly increasing decreasing.  
Are there necessary and sufficient conditions for an arithmetic function $f$ 
to be wildly increasing-decreasing?

Are there sufficient conditions for a sequence  $\mcf = (f_j)_{j=0}^{\infty}$ 
of arithmetic functions to be wildly increasing-decreasing?

\end{document}